\documentclass[a4paper,11pt]{article}

\usepackage{a4wide}
\usepackage{amsthm,amsmath,amssymb}

\def\C{\mathbb C}
\def\eps{\varepsilon}
\def\ga{\gamma}
\def\gb{{\gamma}^{-1}}
\def\GL{{\rm GL}}
\def\GLF{{\rm GL}_2(F)}
\def\GO{{\Gamma}_0({\pi}_F^n)}
\def\GOO{{\Gamma}_0({\pi}_F)}
\def\I2{{\rm I}_2}
\def\ext{{\rm \bf ext}}
\def\MOF{{\rm M}_2({\mathcal O}_F)}
\def\OF{{\mathcal O}_F}
\def\PF{{\pi}_F}
\def\P{{\mathbb P}^1(F)}
\def\Q{\mathbb Q}
\def\sp{\rm Sp}
\def\res{{\rm \bf res}}
\def\v{{\rm val}}

\newtheorem{theo}{Theorem}

\newtheorem{lemma}{Lemma}

\begin{document}

\title{Test vectors for trilinear forms, when  two  representations are unramified and one is special}
\author{Louise Nyssen}
\date{\today}

\maketitle

\section{Introduction}

\bigskip
Let $F$ be a finite extension of  $\Q_p$, with ring of integers $\OF$, and uniformizing parameter $\PF$, whose residual field has  
$q$ elements. For $G=\GLF$, let  $(\pi_1, V_1)$, $(\pi_2, V_2)$ and $(\pi_3, V_3)$ be three irreducible, admissible, infinite dimensional 
representations of  $G$. 
Using the theory of Gelfand pairs, Diprenda Prasad proves in \cite{P} that that the space of $G$-invariant linear forms on $V_1\otimes V_2 \otimes V_3$ 
has dimension at most one. He gives a precise criterion for this dimension to be  one, that we will explain now.

Let $D_F^*$ be the group of invertible elements of  the quaternion division algebra  $D_F$ over $F$. 
When $(\pi_i, V_i)$ is a discrete serie representation of $G$, denote by $(\pi'_i, V'_i)$  the irreducible representation of $D^*_F$ associated 
to $(\pi_i, V_i)$ by the Jacquet-Langlands correspondance. Again, by the theory of Gelfand pairs, 
the space of $D^*_F$-invariant linear forms on $V'_1\otimes V'_2 \otimes V'_3$ has dimension at most one.

Let $\sigma_i$ be the two dimensional representations of the Weil-Deligne group of $F$ associated to the irreducible representations $\pi_i$. 
The triple tensor product $\sigma_1 \otimes \sigma_2 \otimes\sigma_3$ is an eight dimensional symplectic representation of the Weil-Deligne group, 
and has local root number $\eps(\sigma_1 \otimes \sigma_2 \otimes \sigma_3)=\pm 1$. 
When $\eps(\sigma_1 \otimes \sigma_2 \otimes \sigma_3)=- 1$, 
one can prove that the representations $\pi_i$'s are all discrete serie representations of $G$.

\begin{theo} (Prasad, theorem 1.4 of \cite{P} ) Let   $(\pi_1, V_1)$, $(\pi_2, V_2)$, $(\pi_3, V_3)$ be three irreducible, admissible, infinite dimensional 
representations of  $G$ such that the product of their central characters is trivial. 
If all the representations $V_i$'s are cuspidal, 
assume that the residue characteristic of $F$ is not 2. Then \par
$\centerdot$ $\eps(\sigma_1 \otimes \sigma_2 \otimes \sigma_3)=1$ if and only if there exist a non zero $G$-invariant linear form on 
$V_1\otimes V_2 \otimes V_3$ \par
$\centerdot$ $\eps(\sigma_1 \otimes \sigma_2 \otimes \sigma_3)=-1$ if and only if  there exist a non zero $D_k^*$ invariant linear form on 
$V'_1\otimes V'_2 \otimes V'_3$.\par
\end{theo}

Once you got a non zero $G$-invariant linear form $\ell$ on $V_1\otimes V_2 \otimes V_3$, or a non zero  $D_k^*$-invariant linear form $\ell'$ on 
$V'_1\otimes V'_2 \otimes V'_3$, you want to find a vector in  $V_1 \otimes V_2 \otimes V_3$ which is not in the kernel of $\ell$, 
or a vector in  $V'_1\otimes V'_2 \otimes V'_3$ which is not in the kernel of $\ell'$.  Such a vector is called a test vector. 
At first sight, it appears to have strong connections with the new vectors $v_1$, $v_2$ and $v_3$   of the representations $\pi_1$, $\pi_2$ et $\pi_3$.

\begin{theo}\label{vt-000} (Prasad, theorem 1.3 of  \cite{P}) 
When  all the $\pi_i$'s are unramified principal series representations of $G$, $v_1 \otimes v_2 \otimes v_3$ is a test vector for $\ell$. 
\end{theo}

\begin{theo}\label{vt-111} (Gross and Prasad, proposition 6.3 of \cite{GP}) When all   the $\pi_i$'s are  unramified twists of the special representation of $G$ : \par  

$\centerdot$ if $\eps(\sigma_1 \otimes \sigma_2 \otimes \sigma_3)=1$, then $v_1 \otimes v_2 \otimes v_3$ is a test vector for $\ell$,

$\centerdot$ if $\eps(\sigma_1 \otimes \sigma_2 \otimes \sigma_3)=-1$, let $R'$ be the unique maximal order in $D_F$. Then the open 
compact subgroup $R'^* \times R'^*\times R'^*$ fixes a unique line in  $V'_1\otimes V'_2 \otimes V'_3$. Any vector on this line is a test vector for $\ell'$. 
\end{theo} 

The proof by  Gross and Prasad of the first statement of this theorem, actually contains another result:

\begin{theo}\label{vt-110} 
When two of the $\pi_i$'s are  unramified twists of the special representation of $G$ and the third
one belongs to the unramified principal serie of $G$, $v_1 \otimes v_2 \otimes v_3$ is a test vector for $\ell$. 
\end{theo}

But the paper \cite{GP} ends up with  an evidence that $v_1 \otimes v_2 \otimes v_3$ is not always a test vector for $\ell$.
Let $K=\GL (\OF)$ be the maximal compact subgroup of $G$. If  $\pi_1$ and $\pi_2$ are unramified and if $\pi_3$ has conductor $n \geq 1$,
$\ell$ being $G$-invariant, $v_1$ and $v_2$ being $K$-invariant, one gets a $K$-invariant  linear form
$$\left\{ \begin{matrix}
V_3 & \longrightarrow & \C \hfill \cr
v & \longmapsto & \ell(v_1\otimes v_2\otimes v) \cr
\end{matrix}  \right.$$
which must be 0 since $\pi_3$ is ramified. Then $\ell(v_1\otimes v_2\otimes v_3)=0.$

Now Gross and Prasad make the following suggestion. Let $\GO$ be the congruence subgroup
$$\GO = \{  \begin{pmatrix}
a & b \cr
c & d \cr
\end{pmatrix} \in K \quad  \quad c \equiv 0 \quad mod \, {\PF}^n \quad \}$$
and $R$ be a maximal order ${\rm M}_2(F)$ such that $R^* \cap K = \GO$. If $v_2^*$ is a $R^*$-invariant vector in$V_2$, the linear form
$$\left\{ \begin{matrix}
V_3 & \longrightarrow & \C \hfill \cr
v & \longmapsto & \ell(v_1\otimes v_2^* \otimes v) \cr
\end{matrix}  \right.$$
is invariant under the action of $R^* \cap K = \GO$, and one can still hope that $v_1 \otimes v_2^* \otimes v_3$ is  a test vector for $\ell$.
In theorem \ref{vt-001} we will focus on the case $n=1$, and prove that 
$v_1 \otimes v_2^* \otimes v_3$ is  a test vector for $\ell$, up to a condition on $\pi_1$ and $\pi_2$. 
This will almost complete the study of test vectors when the $\pi_i's$ have ramification 0 or 1.

\bigskip
In the long term, the search for test vectors is motivated by the subconvexity problem for $L$-functions.
Roughly speaking, one wants to bound some $L$-functions along the critical line $\Re(z)=\frac{1}{2}$. 
A recent and successful idea in this direction has been to relate triple products of automorphic forms 
to special values of $L$-functions on the critical line. In \cite{BR1} and \cite{BR2} Joseph Bernstein and Andre Reznikov 
did this in the {\it eigenvalue} aspect , and in \cite{V} Akshay Venkatesh did it in the level aspect. 
More details about subconvexity and those related techniques will be found in \cite{MV}. 
Test vectors are key ingredients. Bernstein and Reznikov use an explicit test vector. 
Venkatesh uses a theoretical one, but explains that the bounds would be better with an explicit one (see paragraph 5 of \cite{V}). 
Unfortunately, the difficulty of finding them increases with the ramification of the representations involved.  

There is an extension of Prasad's result in \cite{HS}, where Harris and Scholl prove that 
the dimension of the space of $G$-invariant linear forms on $V_1\otimes V_2 \otimes V_3$ is one when $\pi_1$, $\pi_2$ and $\pi_3$ 
are principal series representations, either irreducible or reducible with their unique irreducible subspace, infinite dimensional. 
They apply the global setting of this to the construction of elements in the motivic cohomology of the product of two modular curves constructed by Beilinson. 

\bigskip
I would like to  thank Philippe Michel for suggesting this problem, 
and Wen-Ching Winnie Li who invited me to spend one semester at PennState University where
I could write the first draft of this paper.

\section{Strategy}

\subsection{Notations}\label{notations}

Let $(\rho, W)$ be a smooth representation of a closed subgroup $H$ of $G$. Let $\Delta_H$ be the modular function on $H$. 
The induction of $\rho$ from $H$ to $G$ is a representation $\pi$ whose space is the space 
${\rm Ind}_{H}^{G}\bigl( \rho \bigr)$  of functions $f$ from $G$ to $W$ 
satisfying the two following conditions :

(1) $\forall h \in H \quad \forall g \in G \quad f(hg)={\Delta_H}^{-\frac{1}{2}}(h) \rho(h) f(g)$,

(2) there exist an open compact subgroup $K_f$ of $G$ such that $$\forall k \in K_f, \quad \forall g \in G, \quad f(gk)= f(g)$$

\noindent where  $G$ acts by right translation. The resulting function will be denoted $\langle \pi(g) , f\rangle$ that is
$$\forall g, g_0 \in G \quad \langle \pi(g),f\rangle(g_0) = f(g_0g).$$ 
With the  additional condition  that $f$ must be compactly supported modulo $H$, 
one gets the {\it compact} induction denoted by ${\rm ind}_{H}^{G}$. 
When $G/H$ is compact, there is no difference between ${\rm Ind}_{H}^{G}$ and ${\rm ind}_{H}^{G}$. 

Let $B$ the Borel subgroup of upper triangular matrices in $G$ and $T$ be the diagonal torus. 
Then we will use $\delta = {\Delta_B}^{-1}$ with 
$\delta \Bigl( \begin{pmatrix} a & b \cr 0 & d \cr \end{pmatrix} \Bigl)  = \vert \frac{a}{d} \vert$ 
and $\Delta_T$ is trivial. The quotient $B \backslash G$ is compact and can be identified with $\P$.

For a smooth representation $V$ of $G$, $V^*$ is the space of linear forms on $V$. 
The contragredient representation $\widetilde{\pi}$ is given by the action of $G$ on  $\widetilde{V}$, the subspace of smooth vectors in $V^*$. 
If $H$ is a subgroup of $G$, $\widetilde{V} \subset \widetilde{V_{\vert H}} \subset V^*$. 

More information about induced and contragredient representations will be found in \cite{BZ}.

\bigskip
Let  $(\pi_1, V_1)$, $(\pi_2, V_2)$ and $(\pi_3, V_3)$ be three irreducible, admissible, infinite dimensional 
representations of  $G$ such that the product of their central characters is trivial. 
Assume that $\pi_1$ and $\pi_2$ are unramified principal series, and that $\pi_3$ has conductor $n\geq 1$. 
Then, according to theorem 1, there exist a non-zero, $G$-invariant linear form $\ell$ on $V_1\otimes V_2 \otimes V_3$, 
and we are looking for a vector $v$ in $V_1\otimes V_2 \otimes V_3$ which is not in the kernel of $\ell$. 
In order to follow Gross and Prasad suggestion, we will consider

$$\ga=
\begin{pmatrix}
\PF^n & 0 \cr
0 & 1 \cr
\end{pmatrix} \qquad {\rm and} \qquad R = {\gb} \MOF \ga .$$

One can easily check that 

$$R^* = \gb K \ga \qquad {\rm and} \qquad R^* \cap K = \GO .$$

If $v_1$, $v_2$ and $v_3$ denote the new vectors of $\pi_1$, $\pi_2$ and $\pi_3$, the vector
$$v_2^*=  \pi_2(\gb) \cdot v_2 $$
is invariant under the action of  $R^*$. Hence we can write 

$$ v_1 \in {V_1}^{K}  \qquad \qquad v_2^*\in {V_1}^{R^*} \qquad \qquad  v_3 \in {V_3}^{R^* \cap K }  $$

According to Gross and Prasad $v_1\otimes v_2^* \otimes v_3$ should be a test vector for $\ell$ , for any $n\geq 1$. 
In this paper, we will focus on the case where $n=1$. We will need the following condition regarding $\pi_1$ and $\pi_2$: 
since they are unramified principal series, they are induced from characters $\chi_1$ and $\chi_2$  of $B$, 
that are required to satisfy

\begin{equation}\label{cond}
\chi_1\Bigl( \begin{matrix}
\PF & 0 \cr
0 & {\PF}^{-1} \cr 
\end{matrix}\Bigr) \not= -1 \quad {\rm or} \quad 
\chi_2\Bigl( \begin{matrix}
\PF & 0 \cr
0 & {\PF}^{-1} \cr 
\end{matrix}\Bigr) \not= -1 
\end{equation}

We will prove 

\begin{theo}\label{vt-001} If $n=1$, and (\ref{cond}) is satisfied, $v_1\otimes v_2^* \otimes v_3$ is a test vector for $\ell$.
\end{theo}

The proof will follow the same pattern as Prasad's proof of theorem \ref{vt-000} in \cite{P}, 
with the necessary changes.

\subsection{Central characters}\label{caractere-central}

Let $\omega_1$, $\omega_2$ and $\omega_3$ be the central caracters of $\pi_1$, $\pi_2$ and $\pi_3$. 
Notice that the condition $\omega_1 \omega_2 \omega_3 =1 $ derives from the $G$-invariance of $\ell$.
Since $\pi_1$ and $\pi_2$ are unramified, $\omega_1$ and $\omega_2$ are unramified too, and so is $\omega_3$
because $ \omega_1 \omega_2 \omega_3 =1 $. 
Let $\eta_i$, for $i\in \{ 1,2,3 \}$ be unramified quasi-characters of $F^*$ with 
$\eta_i^2=\omega_i$ and $\eta_1 \eta_2 \eta _3 =1$. Then
$$V_1 \otimes V_2 \otimes V_3 \simeq 
\bigl(  V_1 \otimes \eta_1^{-1} \bigr) \otimes \bigl(V_2\otimes \eta_2^{-1} \bigr) \otimes \bigl(V_3\otimes \eta_3^{-1} \bigr) $$
as a representation of $G$. Hence it is enough to prove  theorem 4 when the central characters of the representations are trivial.

When $n=1$, it is also enough to prove theorem \ref{vt-001} when $V_3$ is the special representation $\sp$ of $G$ : 
take $\eta_3$ to be the unramified character such that $V_3 = \eta_3 \otimes \sp$.

\subsection{Prasad's exact sequences}\label{suites-exactes}

Let us now explain how Prasad finds $\ell$. It is equivalent to search $\ell$ or to search a non zero element in   
${\rm Hom}_G \Bigl( V_1 \otimes V_2 , \widetilde{V_3} \Bigr)$. 
Since the central characters of $\pi_1$ and $\pi_2$ are trivial, there are unramified characters $\mu_1$ and $\mu_2$ such that  for $i=1$ and $i=2$
 $$\pi_i = {\rm Ind}_{B}^{G} \chi_i \qquad {\rm with} \qquad 
 \chi_i \Bigl( \begin{pmatrix} a & b \cr 0 & d \cr \end{pmatrix} \Bigl) = \mu_i \Bigl(\frac{a}{d} \Bigr)$$

\noindent Hence
$$V_1 \otimes V_2 = {\rm Res}_{G}\,{\rm Ind}_{B \times B}^{G \times G} \Bigl( \chi_1 \times \chi_2 \Bigr)$$
where $G$ is diagonally embedded in $G\times G$ for the restriction. 
The action of $G$ on $B\times B \backslash G\times G = \P \times \P$ has precisely two orbits :  
the first one is  $\{ (u,v) \in \P \times  \P \quad \vert \quad u \not = v  \}$, it is open and  can be identified with $T \backslash G$, 
the second one  is  the diagonal embedding of $\P$ in $\P \times \P$, it is closed and it can be identified with $B \backslash G$. 
Then, we have a short exact sequence of $G$-modules

\begin{equation}\label{courtesuite}
0 \rightarrow {\rm ind}_{T}^{G}\Bigl( \frac{\chi_1}{\chi_2} \Bigr)  \xrightarrow{\ext} V_1 \otimes V_2 \xrightarrow{\res} 
{\rm Ind}_{B}^{G} \Bigl(\chi_1 \chi_2 \delta^{\frac{1}{2}}  \Bigr)\rightarrow 0
\end{equation}

The surjection $\res$ is the restriction of functions from  $G \times G$ to the diagonal part of $B \backslash G \times B \backslash G$, that is
$$\Delta_{B \backslash G} = \Bigl\{ (g,bg) \quad \vert \quad b \in B, \quad g \in G \Bigr\}.$$ 
The injection $\ext$ takes a function 
$f \in {\rm ind}_{T}^{G}\Bigl( \frac{\chi_1}{\chi_2} \Bigr)$ to a function 
$F \in {\rm Ind}_{B \times B}^{G \times G} \Bigl( \chi_1 \times \chi_2 \Bigr)$ given by the relation
$$F \Bigl( g, \begin{pmatrix} 0 & 1 \cr 1 & 0 \cr \end{pmatrix} g  \Bigr) = f(g) \label{rel}.$$
Applying the functor ${\rm Hom}_G \Bigl( \cdot \, , \widetilde{V_3} \Bigr) $, one gets a long exact sequence

\begin{multline}\label{longuesuite} 
0 \rightarrow {\rm Hom}_G \Bigl( {\rm Ind}_{B}^{G} \Bigl(\chi_1 \chi_2 \delta^{\frac{1}{2}}  \Bigr) , \widetilde{V_3} \Bigr) 
  \rightarrow  {\rm Hom}_G \Bigl( V_1 \otimes V_2 , \widetilde{V_3} \Bigr) 
  \rightarrow  {\rm Hom}_G \Bigl( {\rm ind}_{T}^{G}\Bigl( \frac{\chi_1}{\chi_2} \Bigr), \widetilde{V_3} \Bigr) \\
\hfill \downarrow \hskip20mm\\
\hfill \cdots  \leftarrow {\rm Ext}_G^1 \Bigl( {\rm Ind}_{B}^{G} \Bigl(\chi_1 \chi_2 \delta^{\frac{1}{2}}  \Bigr) , \widetilde{V_3} \Bigr)
\end{multline}

\subsection{The simple case}\label{cas-simple}

The situation is easier when $n=1$ and $\mu_1 \mu_2 \vert\cdot\vert^{ \frac{1}{2}} = \vert\cdot\vert^{- \frac{1}{2}}$. 
Then $\pi_3$ is special and there is a natural surjection  
$${\rm Ind}_{B}^{G} \Bigl(\chi_1 \chi_2 \delta^{\frac{1}{2}}  \Bigr) \longrightarrow \widetilde{V_3}$$
whose kernel is the one dimensional subspace of constant functions. Thanks to the exact sequence (\ref{courtesuite}) one gets a surjection 
$$\Psi : V_1 \otimes V_2 \longrightarrow \widetilde{V_3}$$ 
which corresponds to 
$$\ell\left\{ \begin{matrix}
 V_1 \otimes V_2 \otimes V_3 & \longrightarrow & \C \hfill \cr
 v \otimes v' \otimes v'' & \longmapsto & \Psi(v\otimes v') . v'' \cr
\end{matrix}\right.$$ 

The surjection $\Psi$ vanishes on $v_1\otimes v_2^*$ if and only if $\res(v_1\otimes v_2^*)$ has constant value on $\P \simeq B \backslash G$. 
Easy computation proves that it is not constant : the new vectors $v_1$ and $v_2$ are functions from $G$ to $\C$ such that  

$$\forall i \in \{1,2\}, \quad \forall b \in B, \quad \forall k \in K, \quad \quad v_i(bk)= \chi_i(b) \cdot \delta(b)^\frac{1}{2}$$  
and 
$$ \forall g \in G,  \quad \quad v_2^*(g)=v_2(g\gb). $$
Then
$$ ( v_1 \otimes v_2^* ) \Bigl( \begin{pmatrix}
1 & 0 \cr
0 & 1 \cr
\end{pmatrix} \Bigr)  
=  v_1\Bigl( \begin{pmatrix}
1 & 0 \cr
0 & 1 \cr
\end{pmatrix} \Bigr)v_2\Bigl( \gb \Bigr) 
 =  v_2\Bigl(  
\begin{pmatrix}
\PF^{-1} & 0 \cr
0 & 1 \cr
\end{pmatrix}
\Bigr)  =  \mu_2(\PF)^{-1}\vert \PF \vert^{-\frac{1}{2}} =  \frac{\sqrt{q}}{\mu_2(\PF)} $$
and 
$$
( v_1 \otimes v_2^* ) 
\Bigl( \begin{pmatrix}
0 & 1 \cr
1 & 0 \cr
\end{pmatrix} \Bigr)  
=  v_2\Bigl(  \begin{pmatrix}
0 & 1 \cr
1 & 0 \cr
\end{pmatrix}
\begin{pmatrix}
\PF^{-1} & 0 \cr
0 & 1 \cr
\end{pmatrix}
\Bigr)  
= v_2\Bigl(  \begin{pmatrix}
1 & 0 \cr
0 & \PF^{-1} \cr
\end{pmatrix}
\begin{pmatrix}
0 & 1 \cr
1 & 0 \cr
\end{pmatrix}
\Bigr)
=  \frac{\mu_2(\PF)}{\sqrt{q}}. $$

\noindent The representation $\pi_2$ is principal so $\frac{\sqrt{q}}{\mu_2(\PF)} \not = \frac{\mu_2(\PF)}{\sqrt{q}}$ and 
$$ ( v_1 \otimes v_2^* ) \Bigl( \begin{pmatrix}
1 & 0 \cr
0 & 1 \cr
\end{pmatrix} \Bigr)  \not=   ( v_1 \otimes v_2^* ) 
\Bigl( \begin{pmatrix}
0 & 1 \cr
1 & 0 \cr
\end{pmatrix} \Bigr).$$  
Hence, $\Psi$ does not  vanish on $v_1\otimes v_2^*$. 
Then,  $v_1$ being $K$-invariant and $v_2^*$ being $R^*$-invariant, $\Psi(v_1 \otimes v_2^*)$ is a non zero $\GO$-invariant element of $\widetilde{V_3}$, 
that is, a new vector for $\widetilde{\pi_3}$, and it does not vanish on $v_3$ :
$$\ell(v_1\otimes v_2^* \otimes v_3)= \Psi(v_1 \otimes v_2^*).v_3 \not= 0$$
Then $v_1\otimes v_2^* \otimes v_3$ is a test vector for $\ell$.

\subsection{The other case}\label{autre-cas}

If $n \geq 2$ or $\mu_1 \mu_2 \vert\cdot\vert^{ \frac{1}{2}} \not = \vert\cdot\vert^{- \frac{1}{2}}$
then ${\rm Hom}_G \Bigl( {\rm Ind}_{B}^{G} \Bigl(\chi_1 \chi_2 \delta^{\frac{1}{2}}  \Bigr) , \widetilde{V_3} \Bigr) =0$ and by corollary 5.9 of \cite{P} 
$${\rm Ext}_G^1 \Bigl( {\rm Ind}_{B}^{G} \Bigl(\chi_1 \chi_2 \delta^{\frac{1}{2}}  \Bigr) , \widetilde{V_3} \Bigr)=0$$ 
Through the long exact sequence (\ref{longuesuite}) we get an isomorphism 
$${\rm Hom}_G \Bigl( V_1 \otimes V_2 , \widetilde{V_3} \Bigr) \simeq {\rm Hom}_G \Bigl( {\rm ind}_{T}^{G}\Bigl( \frac{\chi_1}{\chi_2} \Bigr), \widetilde{V_3} \Bigr)$$
and by Frobenius reciprocity 
$${\rm Hom}_G \Bigl( {\rm ind}_{T}^{G}\Bigl( \frac{\chi_1}{\chi_2} \Bigr) , \widetilde{V_3} \Bigr) \simeq 
{\rm Hom}_T \Bigl( \Bigl( \frac{\chi_1}{\chi_2} \Bigr) , \widetilde{{V_3}_{\vert T}} \Bigr)$$

\noindent By lemmas 8 and 9 of \cite{W}, this latter space is one dimensional. 
Thus, we have a chain of isomorphic one dimensional vector spaces
$$\begin{matrix}
\ell    & \in & {\rm Hom}_G  \Bigl( V_1 \otimes V_2 \otimes V_3 , \C \Bigr) \\
        &     & \downarrow \wr \\       
\Psi    & \in & {\rm Hom}_G \Bigl( V_1 \otimes V_2 , \widetilde{V_3} \Bigr) \\
        &     & \downarrow \wr \\    
\Phi    & \in & {\rm Hom}_G \Bigl( {\rm ind}_{T}^{G}\Bigl( \frac{\chi_1}{\chi_2} \Bigr), \widetilde{V_3} \Bigr) \\ 
        &     & \downarrow \wr \\    
\varphi & \in & {\rm Hom}_T \Bigl( \Bigl( \frac{\chi_1}{\chi_2} \Bigr) , \widetilde{{V_3}_{\vert T}} \Bigr) \\
\end{matrix}\label{chaine}$$
with generators $\ell$, $\Psi$, $\Phi$ and $\varphi$ corresponding via the isomorphisms. Notice that $\varphi$ is a linear form on $V_3$ such that  
\begin{equation}\label{phi}
\forall t \in T \qquad \forall v \in V_3 \qquad \varphi\bigl(\pi_3(t)v\bigr) = \frac{\chi_2(t)}{\chi_1(t)}\varphi(v) 
\end{equation}

\begin{lemma}\label{lemmeGP}
$\varphi(v_3) \not= 0$.
\end{lemma}

{\it Proof} : this is proposition 2.6 of \cite{GP} with the following translation : \par
- the local field $F$ is the same, \par
- the quadratic extension $K/F$ of Gross and Prasad is $F \times F$ (this case in included in their proof) 
and their group $K^*$ is our torus $T$, \par
- the infinite dimensional representation $V_1$ of  Gross and Prasad is our $\pi_3$,\par
- the one dimensional, unramified representation $V_2$ of  Gross and Prasad is $\frac{\chi_1}{\chi_2}$.\par

Then the representation that Gross and Prasad call $V$ is $\frac{\chi_1}{\chi_2} \otimes \pi_3$ and their condition (1.3) is exactly our condition (\ref{phi}).
The character $\omega$ of Gross and Prasad, which is the central character of their $V_1$, is trivial for us.
Let $\alpha_{K/F}$ be the quadratic character of $F^*$ associated to the extension $K/F$ by local class-field theory,  
and let $\sigma$ and $\sigma_3$ be the  representations of the Weil-Deligne group of $F$ associated to $\frac{\chi_1}{\chi_2}$ and $\pi_3$.   
Thanks to \cite{T} we know that $\eps(\sigma \otimes\sigma_3 ) = \alpha_{K/F}(-1)$   because $K$ is not a field, 
and we are in the first case of proposition 2.6.

\noindent The restriction of $\frac{\chi_1}{\chi_2} \otimes \pi_3$ to the group 
$$M=\Bigl\{  \begin{pmatrix} x & 0 \cr 0 & z \cr \end{pmatrix} \quad \vert  \quad x,y \in \OF^* \Bigr\} \times \GO$$
fixes a unique line in $V_3$ : it is the line generated by the new vector $v_3$. 
According to  Gross and Prasad, a non-zero linear form on $V_3$ which satisfies (\ref{phi}) cannot vanish on $v_3$. 
$\hfill \Box$

\bigskip
\bigskip

We still need to prove that $\ell(v_1\otimes v_2^* \otimes v_3) \not= 0$. For the reason described at the end of section \ref{cas-simple}, 
it is enough to prove that
$$\left\{ \begin{matrix}
 V_3 & \longrightarrow & \C \hfill \cr
 v & \longmapsto &  \ell(v_1\otimes v_2^* \otimes v) \cr
\end{matrix}\right.$$ 
is non zero in $\widetilde{V_3}$. 
\bigskip
In order to do that we want to build a function 
$F$ in $V_1 \otimes V_2$, of the form

\begin{equation}\label{F}
F = \sum_{i\in I} a_i \Bigl< (\pi_1 \otimes \pi_2) (g_i) , v_1 \otimes v_2^*\Bigr> 
\end{equation}

\noindent which vanishes on the closed orbit of $G$ in  $\P \times  \P$. Then, $F$ is in the kernel of $\res$ 
so it is the image by $\ext$ of a function $f \in {\rm ind}_{T}^{G}\Bigl( \frac{\chi_1}{\chi_2} \Bigr)$.
The important point is that $f$ must be the characteristic function of the orbit of the unit 
in the decomposition of $T \backslash G$ under the action of $\GO$, which means : 
\begin{equation}\label{f}
f(g) = \left\{ 
\begin{matrix}
\frac{\chi_1(t)}{\chi_2(t)} & {\rm if} \quad g=tk \quad {\rm with} \quad t \in T  \quad {\rm and} \quad k \in \GO  \cr
0 & {\rm else} \hfill\cr
\end{matrix}
\right.
\end{equation}

Then, the function
$$\left\{ \begin{matrix}
 G & \longrightarrow & \C \hfill \cr
 g & \longmapsto &  f(g) \, \varphi\Bigl( \pi_3(g)v_3 \Bigr)\cr
\end{matrix}\right.   $$
is invariant by the action of $T$ by left translation and  we can do the following computation: on the one hand
\begin{align*} 
\Bigl( \Psi(F) \Bigr) (v_3) & = \Bigl( \Phi(f) \Bigr)(v_3) \\
                            & = \int_{T \backslash G} \! f(g) \, \varphi\Bigl( \pi_3(g)v_3 \Bigr) dg  \\
                            & =\int_{(T\cap K) \backslash\GO} \! \varphi\Bigl( \pi_3(k)v_3 \Bigr) dk \\
                            & = \lambda \cdot \varphi(v_3) .\\
\end{align*}  
where $\lambda$ is a non zero constant. Thanks to lemma \ref{lemmeGP} we know that $\varphi(v_3) \not= 0$ then
$$\Bigl( \Psi(F) \Bigr) (v_3) \not= 0. $$

On the other hand, it comes from (\ref{F}) that 
\begin{align*}
\Bigl( \Psi(F) \Bigr) (v_3) & = \sum_{i\in I} a_i \, \ell \Bigl(\pi_1(g_i)v_1 \otimes \pi_2(g_i) v_2^* \otimes v_3\Bigr)  \\
                            & = \sum_{i\in I} a_i \, \ell \Bigl(v_1 \otimes v_2^* \otimes \pi_3(g_i^{-1})v_3\Bigr) \\
                            & = \Psi(v_1 \otimes v_2^* )\Bigl(\Bigl(\sum_{i\in I} a_i\,  \pi_3(g_i^{-1})\Bigr)v_3\Bigr) \\ 
\end{align*}
then $\Psi(v_1 \otimes v_2^* ) \not= 0$ and  $v_1\otimes v_2^* \otimes v_3$ is a test vector for $\ell$.

\section{Calculations}

\subsection{The big function $F$ and the little function $f$}

The function $F$ has to be $\ext(f)$, where $f$ is the function described by formula (\ref{f}). 
Since $F$ is in $V_1 \otimes V_2 = {\rm Res}_{G}\,{\rm Ind}_{B \times B}^{G \times G} \Bigl( \chi_1 \times \chi_2 \Bigr)$
and $G=BK$, it is enough to know the values of $F$ on  $K \times K$. 

\begin{lemma}\label{lemmeF}$\forall (k,k') \in K \times K$, 

\begin{equation}\label{FK}
F(k,k') = \left\{ 
\begin{matrix}
1 & {\rm if} \quad k \in \GO  \quad {\rm and} \quad  k' \not\in \GOO \cr
0 & {\rm else} \hfill\cr
\end{matrix}
\right.
\end{equation}

\end {lemma}

{\it Proof} :  $F$ must vanish on 
$$\Delta_{B \backslash G} = \Bigl\{ (g,bg) \quad \vert \quad b \in B, \quad g \in G \Bigr\}$$
The other part of $B \backslash G \times B \backslash G$ can be identified with $T \backslash G$ via the bijection 
$$\left\{ \begin{matrix}
\Bigl( B \backslash G \times B \backslash G \Bigr) \setminus \Delta_{B \backslash G} & \rightarrow & T \backslash G \cr
\Bigl( Bg , B\begin{pmatrix} 0 & 1 \cr 1 & 0 \cr \end{pmatrix}g \Bigr) & \longmapsto & Tg \cr
\end{matrix}\right.$$
through which, the orbit of the unit in  $T \backslash G$ under the action of $\GO$ corresponds to 
$$ \Bigl\{  \Bigl( Bk , B\begin{pmatrix} 0 & 1 \cr 1 & 0 \cr \end{pmatrix}k \Bigr) \quad \vert \quad k \in \GO \Bigr\}  $$
Pick any $(k,k') \in K \times K$. If $k' \in Bk$, then 
$k' \in \GO$ if and only if $k \in \GO$, and $k' \in \GOO$ if and only if $k \in \GOO$.
Else, put 
$$k=\begin{pmatrix}a & b \cr c & d\cr \end{pmatrix} \qquad {\rm and} \qquad 
k'=\begin{pmatrix}a' & b' \cr c' & d'\cr \end{pmatrix} .$$
There exist $(b_1, b_2) \in B \times B$ such that 
$$\left\{ \begin{matrix}
k=b_1k_0 \hfill \cr
k'=b_2 \begin{pmatrix}
0 & 1 \cr
1 & 0 \cr
\end{pmatrix}k_0\cr
\end{matrix}\right. \qquad {\rm with} \qquad 
k_0=\begin{pmatrix}c' & d' \cr c & d\cr \end{pmatrix} .$$
Then 
\begin{align*}
k_0 \in \GO  & \iff  c \equiv 0 \quad mod \, {\PF}^n \quad{\rm and}\quad c'd \in \OF^* \cr
             &  \iff c \equiv 0 \quad mod \, {\PF}^n \quad{\rm and}\quad c' \in  \OF^*   \cr
             &   \iff k \in \GO \quad{\rm and}\quad k' \notin \GOO. \cr
\end{align*}
It follows that $(k,k')$ corresponds to an element of the orbit of the unit 
in the decomposition of $T \backslash G$ under the action of $\GO$ if and only if 
$k \in \GO$ and $k' \notin \GOO$.
$\hfill \Box$

\subsection{The big function $F$ when $n=1$}

Now we have to find the coefficients $a_i$ and elements $g_i$ of (\ref{F}) to get the right $F$. 
This can be done for $n=1$. 
For the sake of simplicity, for any family $(g_i)$ of elements of $G$, and $(a_i)$ some complex numbers, denote
$$ \Bigl( \, \sum_{i}a_i \cdot g_i \, \Bigr)( v_1 \otimes v_2^* )  = \sum_{i} a_i \cdot \Bigl< (\pi_1 \times \pi_2 ) (g_i), v_1 \otimes v_2^* \Bigr> $$  
Let $\{ \, \tau_0, \dots \tau_{q-1} \,\}$ be a set of representatives of $\OF / \PF \OF$ in $\OF$, and $A$ be the number
$$A= \Bigl(\frac{\mu_1(\PF)}{\sqrt{q}}- \frac{\sqrt{q}}{\mu_1(\PF)} \Bigr)^{-1} \Bigl(\frac{\mu_2(\PF)}{\sqrt{q}}- \frac{\sqrt{q}}{\mu_2(\PF)} \Bigr)^{-1}$$ 
which can be defined because the representations $\pi_1$ and $\pi_2$ are  principal so 
$\mu_1(\PF)^2 - q \not=0$ and $\mu_2(\PF)^2 - q \not=0$.

\begin{lemma}\label{formule}
When $n=1$ and $1+\mu_1(\PF)^2 \not=0$ the function $F$ is given by

\begin{align*}
F  = \,  A \cdot\Bigl\{ & \quad   \frac{\sqrt{q}}{\mu_2(\PF)}  \cdot \begin{pmatrix}
0 & 1 \cr
\PF & 0 \cr
\end{pmatrix}  +  \frac{\mu_1(\PF)}{\sqrt{q}}  \cdot \begin{pmatrix}
1 & 0 \cr
0 & 1 \cr
\end{pmatrix} \\
&  -  \frac{1}{(1+\mu_1(\PF)^2)} \cdot \frac{\mu_1(\PF)}{\sqrt{q}} \cdot   \frac{\mu_1(\PF)}{\mu_2(\PF)}  \cdot \Bigl(\,  \sum_{i=0}^{q-1}   
\begin{pmatrix}
\PF & {\tau}_i \cr
0 & 1 \cr
\end{pmatrix} +\begin{pmatrix}
0 & 1 \cr
\PF & 0 \cr
\end{pmatrix}\, \Bigr)   \\
&  - \frac{1}{(1+\mu_1(\PF)^2)} \cdot \frac{\mu_1(\PF)}{\sqrt{q}} \cdot  \Bigl( \, \sum_{i=0}^{q-1}   
\begin{pmatrix}
1 & \frac{{\tau}_i}{\PF} \cr
0 & 1 \cr
\end{pmatrix} +\begin{pmatrix}
0 & \frac{1}{\PF} \cr
\PF & 0 \cr
\end{pmatrix}
\, \Bigr)  \quad  \Bigr\}\, ( v_1 \otimes v_2^* )  \\
\end{align*}

When $n=1$ and $1+\mu_2(\PF)^2 \not=0$ the function $F$ is given by

\begin{align*} 
F  = \,   A \cdot \Bigl\{ & \quad    \frac{\sqrt{q}}{\mu_2(\PF)} \begin{pmatrix}
0 & 1 \cr
\PF & 0 \cr
\end{pmatrix}  +  \frac{\mu_1(\PF)}{\sqrt{q}} \begin{pmatrix}
1 & 0 \cr
0 & 1 \cr
\end{pmatrix}  \\
&  -  \frac{1}{(1+\mu_2(\PF)^2)} \cdot \frac{\mu_1(\PF)}{\sqrt{q}}  \cdot \Bigl(\,  \sum_{i=0}^{q-1}   
\begin{pmatrix}
1 & 0\cr
{\tau}_i  & 1 \cr
\end{pmatrix} +\begin{pmatrix}
0 & 1 \cr
1 & 0 \cr
\end{pmatrix}\, \Bigr)   \\
& - \frac{1}{(1+\mu_2(\PF)^2)} \cdot \frac{\mu_2(\PF)}{\sqrt{q}}  \cdot \Bigl( \, \sum_{i=0}^{q-1}   
\begin{pmatrix}
\frac{1}{\PF} & 0 \cr
{\tau}_i & 1 \cr
\end{pmatrix} +\begin{pmatrix}
0 & \frac{1}{\PF} \cr
1 & 0 \cr
\end{pmatrix}
\, \Bigr)  \quad \Bigr\}\, ( v_1 \otimes v_2^* )  \\
\end{align*} 
 
\end{lemma}

{\it Proof} : for $g \in G$ and $k =  \begin{pmatrix}
a & b \cr
c & d \cr
\end{pmatrix}\in K$ in order to compute 
$$\Bigl< \pi_1(g) , v_1 \Bigr> (k) = v_1(kg) \quad {\rm and} \quad  \Bigl< \pi_2(g) , v_2 ^* \Bigr> (k) = v_2(kg\gb)$$
write
$$kg=\begin{pmatrix}
x & y \cr
0 & z \cr
\end{pmatrix} k_1 \quad {\rm and} \quad kg\gb=\begin{pmatrix}
x' & y' \cr
0 & z' \cr
\end{pmatrix} k_2$$
with $k_1$ and $k_2$ in $K$. Then 
$$v_1(kg)= \frac{\mu_1(x)}{\mu_1(z)} \cdot \Bigl\vert\frac{x}{z}\Bigr\vert^{\frac{1}{2}} = \Bigl( \frac{\mu_1(\PF)}{\sqrt{q}}\Bigr)^{(\v\,x -\v\,z)}
\qquad  v_2(kg)=  \Bigl( \frac{\mu_2(\PF)}{\sqrt{q}}\Bigr)^{(\v\,x' -\v\,z')} $$

The following tables give the pairs $ \Bigl( \Bigl< \pi_1(g) , v_1 \Bigr> (k),  \Bigl< \pi_2(g) , v_2 ^* \Bigr> (k) \Bigr)$. The entries, are : 
an element $g$ in $G$, $val(c)$ and $val(d)$ where $(c,d)$ is the second line of $k$.

The first table is inspired by the formula
$$T_{\PF}= 
K \begin{pmatrix}
\PF & 0 \cr
0 & 1 \cr
\end{pmatrix}K = \sqcup_{i=1}^{q-1} 
\begin{pmatrix}
\PF & \tau_{i} \cr
0 & 1 \cr
\end{pmatrix}K + 
\begin{pmatrix}
0 & 1 \cr
\PF & 0 \cr
\end{pmatrix}K$$
$$\begin{array}{|l|c|c|}
\hline &&\cr
\hfill g  \hfill & \v(c) = 0 & \v(c) \geq 1  \cr
&&\cr \hline &&\cr
 \begin{pmatrix}
\PF & \tau_{i} \cr
0 & 1 \cr
\end{pmatrix} \quad  {\rm such} \quad {\rm that} \quad  c \tau_{i}+d \in {\OF}^*  & \Big(\frac{\mu_1(\PF)}{\sqrt{q}}\Big),1 & \Big(\frac{\mu_1(\PF)}{\sqrt{q}}\Big),1  \cr
&&\cr \hline &&\cr
\begin{pmatrix}
\PF & \tau_{i_0} \cr
0 & 1 \cr
\end{pmatrix} \quad  {\rm such} \quad {\rm that} \quad c \tau_{i_0}+d \in \PF \OF     & \Big(\frac{\mu_1(\PF)}{\sqrt{q}}\Big)^{-1},1  & \emptyset \cr
&&\cr \hline &&\cr
\begin{pmatrix}
0 & 1 \cr
\PF & 0 \cr
\end{pmatrix}  & \Big(\frac{\mu_1(\PF)}{\sqrt{q}}\Big),1    & 
\Big(\frac{\mu_1(\PF)}{\sqrt{q}}\Big)^{-1},1  \cr
&&\cr \hline
\end{array}$$

\noindent
Fix $$F_1 = \Bigl(\,  \sum_{i=0}^{q-1}   
\begin{pmatrix}
\PF & {\tau}_i\cr
0  & 1 \cr
\end{pmatrix} +\begin{pmatrix}
0 & 1 \cr
\PF & 0 \cr
\end{pmatrix}\, \Bigr)(v_1 \otimes v_2^*)$$
It comes out that $\forall (k,k') \in K \times K$
$$  F_1 (k,k')= q . \frac{\mu_1(\PF)}{\sqrt{q}} + \frac{\sqrt{q}}{\mu_1(\PF)} = 
\frac{\sqrt{q}}{\mu_1(\PF)} \cdot (1+\mu_1(\PF)^2)$$
Now consider $\gb F_1 (k,k')$. On the one hand 

\begin{align*}  
\gb F_1 & = \Bigl(\,  \sum_{i=0}^{q-1}   
\gb \begin{pmatrix}
\PF & {\tau}_i\cr
0  & 1 \cr
\end{pmatrix} + \gb \begin{pmatrix}
0 & 1 \cr
\PF & 0 \cr
\end{pmatrix}\, \Bigr)(v_1 \otimes v_2^*) \\
& =  \Bigl( \, \sum_{i=0}^{q-1}   
\begin{pmatrix}
1 & \frac{{\tau}_i}{\PF}\cr
0 & 1 \cr
\end{pmatrix} +\begin{pmatrix}
0 & \frac{1}{\PF} \cr
\PF & 0 \cr
\end{pmatrix}
\, \Bigr) (v_1 \otimes v_2^*).
\end{align*}

\noindent 
On the other hand, for any $(k,k')$ in $K \times K$, $$\Bigl( \gb F_1 \Bigr) (k,k')=F_1(k\gb,k'\gb).$$
Take $k=\begin{pmatrix}
a & b \cr
c & d \cr
\end{pmatrix}$, and $ k \gb= \begin{pmatrix}
\frac{a}{\PF} & b \cr
\frac{c}{\PF} & d \cr
\end{pmatrix}$. If $\v\,c=0$, then $\v\,d + 1 \geq \v\,c $ and
$$ k \gb =  \begin{pmatrix}
\frac{ad-bc}{c} & \frac{a}{\PF} \cr
0 &  \frac{c}{\PF} \cr
\end{pmatrix}\begin{pmatrix}
0 & -1 \cr
1 & \frac{\PF d}{c} \cr
\end{pmatrix}$$
with $$ \begin{pmatrix}
0 & -1 \cr
1 & \frac{\PF d}{c} \cr
\end{pmatrix} \in K \quad {\rm and} \quad (\chi_1 \cdot {\delta}^\frac{1}{2})\begin{pmatrix}
\frac{ad-bc}{c} & \frac{a}{\PF} \cr
0 & \frac{c}{\PF} \cr
\end{pmatrix} = \frac{\mu_1(\PF)}{\sqrt{q}}$$ 
If $\v\,c \geq 1$, then $\v\,d = 0 $ and
$$ k \gb =  \begin{pmatrix}
\frac{ad-bc}{\PF d} & b \cr
0 &  d \cr
\end{pmatrix}\begin{pmatrix}
1 & 0 \cr
\frac{c}{\PF d} &  1 \cr
\end{pmatrix}$$
with $$ \begin{pmatrix}
1 & 0 \cr
\frac{c}{\PF d} &  1 \cr
\end{pmatrix} \in K \quad {\rm and} \quad (\chi_1 \cdot {\delta}^\frac{1}{2})(\begin{pmatrix}
\frac{ad-bc}{\PF d} & b \cr
0 &  d \cr
\end{pmatrix}) = \frac{\sqrt{q}}{\mu_1(\PF)}$$ 
The same calculation with $k'$ leads to the following:

\begin{equation*}
F_1( k\gb ,k'\gb ) = \frac{\sqrt{q}}{\mu_1(\PF)} \cdot (1+\mu_1(\PF)^2) \cdot \left\{ 
\begin{matrix}
\frac{\mu_1(\PF)}{\sqrt{q}}\cdot\frac{\mu_2(\PF)}{\sqrt{q}} & {\rm if} \quad  \v\,c=\v\,c'=0 \hfill \cr
\frac{\mu_1(\PF)}{\sqrt{q}}\cdot\frac{\sqrt{q}}{\mu_2(\PF)} & {\rm if} \quad  \v\,c=0 \quad {\rm and} \quad \v\,c'\geq 1 \hfill \cr
\frac{\sqrt{q}}{\mu_1(\PF)}\cdot\frac{\mu_2(\PF)}{\sqrt{q}} & {\rm if} \quad  \v\,c \geq 1 \quad {\rm and} \quad \v\,c'=0 \hfill \cr
\frac{\sqrt{q}}{\mu_1(\PF)}\cdot\frac{\sqrt{q}}{\mu_2(\PF)} & {\rm if} \quad  \v\,c \geq 1 \quad {\rm and} \quad \v\,c'\geq 1 \hfill \cr
\end{matrix}
\right.
\end{equation*}

Now, with  the simple table
$$\begin{array}{|l|c|c|}
\hline &&\cr
\hfill g \hfill & \v(c) = 0  & \v(c) \geq 1  \cr
&&\cr \hline &&\cr
 \begin{pmatrix}
1 & 0 \cr
0 & 1 \cr
\end{pmatrix}  & 1,\Big(\frac{\mu_2(\PF)}{\sqrt{q}}\Big)    & 
1,\Big(\frac{\mu_2(\PF)}{\sqrt{q}}\Big)^{-1} \cr
&&\cr \hline &&\cr
\begin{pmatrix}
0 & 1 \cr
\PF & 0 \cr
\end{pmatrix}  & \Big(\frac{\mu_1(\PF)}{\sqrt{q}}\Big),1    & 
\Big(\frac{\mu_1(\PF)}{\sqrt{q}}\Big)^{-1},1  \cr
&&\cr \hline 
\end{array}$$
and 
\begin{align*}  
F  = &  A \cdot \Bigl( \frac{\sqrt{q}}{\mu_2(\PF)} \cdot \begin{pmatrix}
0 & 1 \cr
\PF & 0 \cr
\end{pmatrix}  +  \frac{\mu_1(\PF)}{\sqrt{q}}  \cdot \begin{pmatrix}
1 & 0 \cr
0 & 1 \cr
\end{pmatrix} \Bigr) (v_1 \otimes v_2^*) \\
& -\frac{A}{(1+\mu_1(\PF)^2)} \cdot \frac{\mu_1(\PF)}{\sqrt{q}} \cdot \Bigl( \frac{\mu_1(\PF)}{\mu_2(\PF)}  \cdot F_1 + \gb F_1 \Bigr)  \\
\end{align*}   \par
one gets the first formula of lemma \ref{formule}.

\bigskip
The second formula of lemma \ref{formule} is obtained by considering the  decomposition
$$K = \sqcup_{j=0}^{q-1} 
\begin{pmatrix}
1 & 0 \cr
\tau_{j} & 1 \cr
\end{pmatrix} \GOO \sqcup 
\begin{pmatrix}
0 & 1 \cr
1 & 0 \cr
\end{pmatrix}\GOO$$ 
Then, using the table
$$\begin{array}{|l|c|c|c|}
\hline &&& \cr
                &  \v(c) = 0    &                   &  \v(c) \geq 1     \cr
\hfill g \hfill &               & \v(d) = \v(c) = 0 &                   \cr
                &  \v(d) \geq 1 &                   &  \quad \v(d)  = 0 \cr
&&& \cr \hline &&& \cr
\begin{pmatrix}
1 & 0 \cr
0 & 1 \cr
\end{pmatrix}  
& 1,\Big(\frac{\mu_2(\PF)}{\sqrt{q}}\Big) & 1,\Big(\frac{\mu_2(\PF)}{\sqrt{q}}\Big) & 1,\Big(\frac{\mu_2(\PF)}{\sqrt{q}}\Big)^{-1}  \cr
&&& \cr \hline &&& \cr
 \begin{pmatrix}
       1 & 0 \cr
\tau_{j} & 1 \cr
\end{pmatrix}  \quad 
\begin{array}{l}  {\rm such} \quad  {\rm that} \quad \tau_j \not= 0 \cr {\rm and} \quad   d\tau_{j}+c \in {\OF}^*  \cr \end{array}
  & 1,\Big(\frac{\mu_2(\PF)}{\sqrt{q}}\Big)    & 1,\Big(\frac{\mu_2(\PF)}{\sqrt{q}}\Big)  & 
1,\Big(\frac{\mu_2(\PF)}{\sqrt{q}}\Big)  \cr
&&& \cr \hline &&& \cr
\begin{pmatrix}
1 & 0 \cr
\tau_{j_0} & 1 \cr 
\end{pmatrix} 
\begin{array}{l}  {\rm such} \quad  {\rm that} \quad \tau_j \not= 0 \cr {\rm and} \quad d\tau_{j_0}+c \in \PF \OF  \cr \end{array}
& \emptyset  & 1,\Big(\frac{\mu_2(\PF)}{\sqrt{q}}\Big)^{-1}  & \emptyset \cr
&&& \cr \hline& && \cr
\begin{pmatrix}
0 & 1 \cr
1 & 0 \cr
\end{pmatrix} & 
1, \Big(\frac{\mu_2(\PF)}{\sqrt{q}}\Big)^{-1} & 1,\Big(\frac{\mu_2(\PF)}{\sqrt{q}}\Big) 
& 1, \Big(\frac{\mu_2(\PF)}{\sqrt{q}}\Big)  \cr
&&& \cr \hline
\end{array}$$

\noindent
one gets a function 
$$F_2 = \Bigl(\,  \sum_{i=0}^{q-1}   
\begin{pmatrix}
1 & 0\cr
{\tau}_i  & 1 \cr
\end{pmatrix} +\begin{pmatrix}
0 & 1 \cr
1 & 0 \cr
\end{pmatrix}\, \Bigr) (v_1 \otimes v_2^*)$$
which satisfies $\forall (k,k') \in K \times K$
$$  F_2 (k,k')= q . \frac{\mu_2(\PF)}{\sqrt{q}} + \frac{\sqrt{q}}{\mu_2(\PF)} = \frac{\sqrt{q}}{\mu_2(\PF)} \cdot (1+\mu_2(\PF)^2).$$
This is the same situation as the previous one : by computing $\gb F_2$ and choosing the right coefficients, one gets the second formula of lemma \ref{formule}.
$\hfill \Box$
 
\vskip2cm
\noindent {\bf {\Large Conclusion :}}
Thus,  we could write the function $F$ for $n=1$ and $1+\mu_1(\PF)^2 \not= 0$ or  $1+\mu_2(\PF)^2 \not= 0$. 
The latter condition is precisely condition \ref{cond} of theorem \ref{vt-001}, which is now proved.  
Of course, it would be interesting  to remove this condition and then to find $F$ for any $n$.


\vskip2cm




\begin{thebibliography}{1}

\bibitem[B-Z]{BZ} Joseph Bernstein and Andrei Zelevinsky, {\it Representations of the group GL(n,F) where F is a non-archimedian local field}.  
Russian Mathematical Surveys {\bf 31:3} (1976), 1-68.

\bibitem[B-R 1]{BR1} Joseph Bernstein and Andre Reznikov, {\it Estimates of automorphic functions}.  Moscow Mathematic Journal  {\bf 4}, no.1 (2004), 
19-37.

\bibitem[B-R 2]{BR2} Joseph Bernstein and Andre Reznikov, {\it Periods, subconvexity and representation theory}.  Journal of differential
geometry {\bf 70} (2005), 129-142.

\bibitem[G-P]{GP} Benedict H.Gross and Diprenda Prasad, {\it Test Vectors for Linear forms}. Mathematische Annalen {\bf 291} (1991), 343-355.

\bibitem[H-S]{HS} Michael Harris and Anthony Scholl, {\it A note on trilinear forms for reducible representations and Beilinson conjectures}.  Journal of the European
Mathmatical Society {\bf 2001, 1} (2001), 93-104.


\bibitem[M-V]{MV} Philippe Michel and Akshay Venkatesh, {\it Equidistribution, $L$-functions and Ergodic theory : on some problem of Yu. V.
Linnik}.  Preprint {\bf } (2005).

\bibitem[P]{P} Diprenda Prasad, {\it Trilinear forms for representations of GL(2) and local $\eps$-factors}. 
Composotio Mathematica {\bf 75} (1990), 1-46.

\bibitem[T]{T} J. Tunnell, {\it Local $\eps$-factors and characters of GL(2)}.  American Journal of Mathematics {\bf 105} (1983), 1277-1308.

\bibitem[V]{V} Akshay Venkatesh, {\it Sparse equidistribution problems, period bounds, and subconvexity}. Preprint {\bf } (2005).

\bibitem[W]{W} Jean-Loup Waldspurger, {\it Sur les valeurs de certaines fonctions $L$ automorphes en leur centre de sym\'etrie}.  
Compositio Mathematica {\bf 54} (1985), 173-242.

\end{thebibliography}
\end{document}